\documentclass[twoside,11pt,leqno]{article}
\usepackage{amssymb}
\usepackage{amsthm}
\usepackage{graphicx} 
\usepackage{mathrsfs} 
\usepackage{multirow}  
\usepackage{amsmath}

\newcommand{\Q}{\mathbb Q}
\newcommand{\RR}{\mathbb R}
\newcommand{\NN}{\mathbb N}
\newcommand{\CC}{\mathbb C}

\def\F{{\mathcal F}}

\def\U{{\mathcal{U}}}
\def\Sym{{\rm Sym}}

\def\0{{\emptyset}}

\newtheorem{thm}{Teorema}[section]

\newtheorem{prop}[thm]{Proposici\'on}

\theoremstyle{definition}

\newtheorem{ejem}[thm]{Ejemplo}
\newtheorem{defi}[thm]{Definici\'on}

\newcommand{\rnf}{\renewcommand{\thefootnote}{\arabic{footnote}}}

\title{Aproximaci\'on m\'etrica de grupos: una breve perspectiva\footnote{Este art\'iculo es parte de la tesis de licenciatura de la  segunda autora bajo la direcci\'on del primer autor. La tesis se present\'o en la Unidad Acad\'emica de Matem\'aticas de la Universidad Aut\'onoma de Zacatecas en marzo del 2016 como requisito para obtener el t\'itulo de Licenciada en Matem\'aticas.}
}

\author{\rnf   Luis Manuel Rivera
\and  \rnf     Nidya Monserrath Veyna Garc\'ia
}

\date{}

\begin{document}
\maketitle






\begin{abstract} \noindent
Los grupos s\'oficos y los grupos hiperlineales han generado una gran cantidad de investigaci\'on en los \'ultimos a\~nos en diversas \'areas de matem\'aticas tales como teor\'ia geom\'etrica de grupos, din\'amica simb\'olica y \'algebra de operadores. Adem\'as, los grupos s\'oficos han ganado inter\'es porque se ha  demostrado que cumplen varias conjeturas a\'un abiertas para los grupos en general. Las definiciones de ambas clases de grupos son an\'alogas y se pueden pensar como dentro de una clase de grupos de reciente estudio que se conocen como los grupos que tienen la propiedad de apro\-xi\-ma\-ci\'on m\'etrica. En este art\'iculo se presenta un panorama general de dicha clase de grupos.
\end{abstract}

\noindent 


\section{Introducci\'on}

El concepto de ultraproducto aparece de manera general en un art\'iculo de Jerzy \L o\v{s} en 1955 \cite{jerzy}, pero fue en a\~nos recientes cuando ha tomado importancia en ciertas \'areas de la teor\'ia de grupos. Lo anterior porque se ha demostrado que varias clases  importantes de grupos tienen la ca\-rac\-te\-r\'is\-ti\-ca com\'un de ser isomorfos a subgrupos de ciertos ultraproductos de grupos, y en este caso decimos que el grupo en cuesti\'on se aproxima por los grupos con los cuales se construye el ultraproducto. Fue en los a\~nos noventa cuando se comenzaron a definir estas clases de grupos, siendo los grupos s\'oficos y los grupos hiperlineales dos de las clases m\'as importantes, las cuales surgen en diferentes ramas de la matem\'atica (din\'amica simb\'olica y \'algebra de operadores, respectivamente). Estos grupos han ganado inter\'es porque se ha demostrado que cumplen varias conjeturas importantes que  siguen abiertas para grupos en general. 

La clase de grupos s\'oficos fue definida, sin usar ultraproductos, por Gromov \cite{gromov-artc}, quien demostr\'o que esta clase de grupos satisfacen la conjetura de sobreyuntividad de Gottschalk \cite{gott}. El nombre de grupos s\'oficos se debe a Weiss \cite{weiss}. Estos grupos son una generalizaci\'on com\'un de los grupos residualmente finitos y de los grupos amenables, y tambi\'en incluyen, entre otros, a  los grupos encajables localmente en finito, estos \'ultimos presentados por Vershik y Gordon \cite{VeG}.

La clase de grupos hiperlineales fue definida por R\u{a}dulescu \cite{radu} quien demostr\'o que estos grupos cumplen la conjetura del encaje de Connes en su versi\'on para grupos. Elek y Szab\'o \cite{ES2} demostraron que los grupos s\'oficos son hiperlineales al mostrar que se pueden encajar en un cierto ultraproducto m\'etrico de grupos sim\'etricos finitos, lo que implica que los  grupos s\'oficos cumplen la conjetura del encaje de Connes. 

Elek y Szab\'o \cite{ES1} tambi\'en demuestran que los grupos s\'oficos cumplen con la conjetura de finitud directa de Kaplansky y con la conjetura del determinante de L\"uck \cite{ES2}. Posteriormente, Thom \cite{thom-conj} demostr\'o que dichos grupos tamb\'ien cumplen con la conjetura de autovalores algebraicos de J. Dodziuk, P. Linnell, V. Mathai, T. Schick y S. Yates. A la fecha no se conocen ejemplos de grupos que no  sean s\'oficos o hiperlineales. M\'as detalles sobre estas dos clases de grupos se pueden consultar, por ejemplo, en el resumen de Pestov \cite{pes}, en el libro de Capraro y Lupini \cite{capraro-lup} y en la monograf\'ia de Ceccherini-Silberstein y Coornaert \cite{b1}.

 A partir de entonces, se han definido otras clases de grupos como posibles generalizaciones de los grupos s\'oficos: los grupos s\'oficos d\'ebiles, los grupos s\'oficos lineales y los grupos $K$-lineales. El siguiente diagrama muestra la relaci\'on entre algunas de estas clases de grupos:\

\begin{center}
\scalebox{0.8}{
\begin{tabular}{ccccccccc}
Finito       & $\Rightarrow$ & Res. Finito   & $\Rightarrow$ &     LEF      &                   & S\'ofico Lineal & $\Rightarrow$ & S\'ofico D\'ebil\\
$\Downarrow$ &               & $\Downarrow$  &               & $\Downarrow$ &                   &   $\Uparrow$  &               &  \\
Amenable     & $\Rightarrow$ & Res. Amenable & $\Rightarrow$ &      LEA     & $\Longrightarrow$ & S\'ofico        & $\Rightarrow$ & Hiperlineal \\
\end{tabular}
}
\end{center}

Los grupos s\'oficos, los grupos s\'oficos d\'ebiles, los grupos s\'oficos lineales y los grupos $K$-lineales se pueden definir sin usar ultraproductos, mediante definiciones que tiene muchas similitudes. La diferencia principal entre las respectivas  definiciones para cada tipo de grupos, radica en la clase de grupos m\'etricos que aproximan al grupo en cuesti\'on. Por ejemplo, los grupos s\'oficos se aproximan por grupos sim\'etricos finitos equipados con la m\'etrica de Hamming y los grupos hiperlineales se aproximan por la clase de grupos unitarios de rango finito equipados con la m\'etrica de Hilbert-Schmidt.

En este art\'iculo se presenta un breve panorama de esta clase de grupos. Nuestro tratamiento del tema es limitado e incompleto pero tiene como objetivo presentar el \'area de investigaci\'on sobre la  aproximaci\'on de grupos de manera breve. A pesar de que la literatura sobre este tema est\'a en constante crecimiento y se cuenta con al menos tres fuentes en donde se resumen algunos de los resultados del \'area (\cite{b1, capraro-lup, pes}),  en la actualidad, a conocimiento de los autores, no se cuenta con literatura en espa\~nol sobre este tema. Adem\'as, a diferencia de \cite{b1, capraro-lup, pes},  en este art\'iculo se aborda el tema de aproximaci\'on de manera general como se ha hecho recientemente por varios autores \cite{lev-gleb, hayes2, holt-rees, artstolz, thom2, stolz-thom}. La mayor\'ia de resultados se presentan sin demostraci\'on pero se indica las fuentes en donde se pueden consultar. 

El resto del art\'iculo est\'a organizado como sigue. En la secci\'on \ref{sec-ultrafiltros} se presenta un breve resumen de definiciones y resultados sobre ultrafiltros. En la secci\'on \ref{sec-ultra-grupos} se presentan a los ultraproductos y se define a los grupos que son localmente encajables en una clase de grupos por medio de dos definiciones, una de las cuales usa ultraproductos. En la secci\'on \ref{sec-aproximacion-grupos} se define a los ultraproductos m\'etricos y se definen el concepto de aproximaci\'on m\'etrica de grupos. Posteriormente se dan varios ejemplos de esta clase de grupos: los grupos s\'oficos, los grupos hiperlineales, los grupos s\'oficos d\'ebiles y los grupos $K$-lineales. En la secci\'on \ref{grupos-soficos} se extiende la exposici\'on sobre los grupos s\'oficos, que es la clase m\'as estudiada en la teor\'ia de aproximaci\'on de grupos. Se presentan varias definiciones equivalentes de grupos s\'oficos, algunos ejemplos y algunas propiedades de cerradura de esta clase de grupos. En la secci\'on \ref{2-prop-grupos-aprox} se presenta una segunda definici\'on de grupos aproximables que no usa ultrafiltros, esta definici\'on es an\'aloga a la de grupos s\'oficos. Adem\'as se presentan algunos resultados de cerradura de esta clase de grupos.
 
\section{Ultrafiltros}\label{sec-ultrafiltros}
Vamos a presentar una breve introducci\'on a la teor\'ia de ultrafiltros. Esta secci\'on esta basada principalmente en \cite[cap\'itulo 2]{ak-khan}, \cite[ap\'endice J]{b1} y \cite[cap\'itulo 1]{b2}.\\
\begin{defi}
Sea $X$ un conjunto no vac\'io. Un {\it filtro} $\mathcal{F}$ en $X$ es una familia no vac\'ia de subconjuntos de $X$ que satisface
\begin{enumerate}
\item [{\it F1})] \label{con1fil} $\emptyset \notin \mathcal{F}$;
\item [{\it F2})]\label{con2fil} Si $A,B \in \mathcal{F}$ entonces $A \cap B \in \mathcal{F}$;
\item [{\it F3})] \label{con3fil} Si $A \in \mathcal{F}$ y $A \subseteq B$ entonces $B \in \mathcal{F}$.
\end{enumerate}
\end{defi}

Algunos autores, por ejemplo en \cite{ak-khan}, no incluyen la condici\'on~(F1) en la definici\'on de filtro. Notemos que si usamos \'unicamente las condiciones (F2) y (F3) en la definici\'on anterior, sigue que $\emptyset \in \mathcal{F}$ si y solo si $\mathcal{F} = \mathcal{P}(X)$. En este caso, al conjunto $\mathcal{P}(X)$ se le llama {\it filtro impropio}. 

Los siguientes son ejemplos de filtros:
\begin{ejem}
\begin{enumerate}
\item[(a)] Para $x_0 \in X$, la familia $\F_{x_0}=\{A \subseteq X \colon x_0 \in A\}$ es un filtro en $X$ y se llama {\it filtro principal} generado por $x_0$.
\item[(b)] Sea $X$ un conjunto infinito. La familia \[\F =\{A \subseteq X \colon X \setminus A \text{ es finito}\},\] es un filtro en $X$ y se llama el {\it filtro de Fr\'echet} en $X$. Notemos que este filtro no es principal.
\item[(c)] Sea $T$ una topolog\'ia en $X$, y $x \in X$. El conjunto 
\[
\mathcal{N}(x)= \{V \colon V \text{ es una vecindad de }  x \}
\]
es un filtro en $X$. A este filtro se le conoce como el {\it filtro de vecindades} de $x$.
\item[(d)] Sea $(X, \leq)$ un conjunto dirigido no vac\'io. Un subconjunto $A \subseteq X$ se llama {\it residual} en $X$ si existe $x_0 \in X$ tal que el conjunto $\{x \in X: x_0 \leq x\}$ es un subconjunto de $A$. El conjunto $\F_r(X)$ de todos los subconjuntos residuales de $X$ es un filtro en $X$ y se llama el {\it filtro residual} en $X$.

\end{enumerate}
\end{ejem}

Notemos que si consideramos al conjunto dirigido $(\mathbb{N}, \leq)$, un subconjunto $A \subseteq \mathbb{N}$ es residual si y solo si $\mathbb{N} \setminus A$ es finito. Entonces, se tiene que el filtro residual en $\mathbb{N}$ es el filtro de Fr\'echet en $\mathbb{N}$. 

\begin{prop}
\label{4}
Un ultrafiltro $\U$ en $X$ es principal si y solo si existe un conjunto finito en $\U$.
\end{prop}

\begin{defi}
Sea $X$ un conjunto. Decimos que una familia $\mathcal{A}$ de subconjuntos de $X$ tiene la {\it propiedad de intersecci\'on finita}, o que es un {\it sistema centrado}, si la intersecci\'on de cualquier subcolecci\'on finita de $\mathcal{A}$ es no vac\'ia. 
\end{defi}
\begin{prop}
\label{PIF}
Sea $X$ un conjunto no vac\'io y sea $\mathcal{A} \subset \mathcal{P}(X)$. Entonces, existe un filtro en $X$ que contiene a $\mathcal{A}$ si y solo si $\mathcal{A}$ tiene la propiedad de intersecci\'on finita.
\end{prop}


\begin{defi}
\label{teoultra}
Un filtro $\mathscr{F}$ en $X$ se llama {\it ultrafiltro} si es un filtro maximal (relativo a la inclusi\'on).
 \end{defi}

\begin{prop}
\label{2}
Un filtro $\U$ en $X$ es un ultrafiltro si y solo si para todo $A \subseteq X$, o bien $A \in \U$, o bien $X \setminus A \in \U$.
\end{prop}
\begin{defi}
Un ultrafiltro $\U$ es libre, si sus elementos no tienen puntos en com\'un, es decir, si  $\bigcap \U = \emptyset$.
\end{defi}
Observemos que cualquier filtro principal es un ultrafiltro. Si consideramos a $X = \mathbb{N}$ tenemos que el filtro de Fr\'echet en $\mathbb{N}$ no es un ultrafiltro porque ninguno de los conjuntos  $2\mathbb{N}$ o $\mathbb{N} \setminus 2\mathbb{N}$ pertenecen a dicho filtro, sin embargo se sabe que cualquier ultrafiltro libre sobre $X$ contiene al filtro de Fr\'echet. 
La siguiente proposici\'on requiere del Lema de Zorn. 
\begin{thm}\label{teo-ultra}
Sea $X$ un conjunto no vac\'io. Cualquier filtro en $X$ es un subconjunto de alg\'un ultrafiltro en $X$.
\end{thm}

\begin{prop}
\label{prinlibre}
Todo ultrafiltro es, o bien principal, o bien libre.
\end{prop}

La existencia de los ultrafiltros libres es una consecuencia del Lema de Zorn. Por ejemplo, 
considere la siguiente colecci\'on de subconjuntos de $\mathbb{N}_+$ 
\[
\mathcal{C} = \{\mathbb{N}_+ \setminus \{n\}: n \in \mathbb{N}_+ \}.
\]
Esta colecci\'on de conjuntos tiene la propiedad de intersecci\'on finita, y es tal que $\bigcap \mathcal{C} = \0$. As\'i por la proposici\'on \ref{PIF} y el teorema \ref{teo-ultra}, existe un ultrafiltro $\U$ que lo contiene. Adem\'as, $\cap \U \subseteq \cap \mathcal{C}=\emptyset$. Por lo tanto $\U$ es un ultrafiltro libre.

\begin{defi}
Sea $X$ un conjunto, $Y$ un espacio topol\'ogico, $y_0 \in Y$, $f:X\rightarrow Y$ una funci\'on y $\F$ un filtro en $X$. Se dice que $y_0$ es un {\it l\'imite} de $f$ a lo largo de $\F$ (o que $f(x)$ {\it converge} a $y_0$ a lo largo de $\F$),  si $f^{-1}(V)=\{x \in X\colon f(x) \in V\}$ pertenece a $\F$ para todo vecindad $V$ de $y_0$.  Si adem\'as $\F$ es un ultrafiltro decimos que  $y_0$ es un {\it ultral\'imite} de $f$ a lo largo de $\F$. Si $y_0$ es un l\'imite \'unico se denota por  
 \[
 \lim_{x \to \F}f(x)=y_0, \; \text{ o simplemente por  }  
 \lim_{\F}f(x)=y_0.
 \]

\end{defi}

\begin{prop}
Sean $X$ un conjunto, $Y$ un espacio topol\'ogico compacto, $f \colon X \to Y$ una funci\'on, y $\U$ un ultrafiltro en $X$. Entonces existe $y_0 \in Y$ tal que $f(x)$ converge a $y_0$ a lo largo de $\F$. Adem\'as, si $Y$ es Hausdorff, entonces $y_0$ es \'unico.
\end{prop}

La definici\'on de ultral\'imite para el caso de series de n\'umeros reales, queda como sigue:
\begin{defi}
Sea $\{x_i\}_{i \in I}$ una familia de n\'umeros reales y $\;\U$ un ultrafiltro en $I$. Decimos que $ \lim\limits_{u} x_i = x \in \mathbb{R}$ si para todo $\epsilon > 0$ se tiene $$\{i\in I: |x_i - x| < \epsilon\} \in \U$$
\end{defi}
El c\'alculo del ultral\'imite depende de la selecci\'on del ultrafiltro. Por ejemplo, 
sea $\{x_n\}_{n \in \mathbb{N}}$ una secuencia convergente de n\'umeros reales. Sea $\U$ un ultrafiltro principal en $\mathbb{N}$ generado por $\{n_0\}$, se puede demostrar que $\lim\limits_{u} x_n = x_{n_0}$.

Ahora, sea $\U$ un ultrafiltro libre en $\NN$, sea $x$ el l\'imite cl\'asico de $\{x_n\}_{n \in \mathbb{N}}$, se puede demostrar que $ \lim\limits_{u} x_n = x$. 


\section{Ultraproductos de grupos}\label{sec-ultra-grupos}
Vamos a presentar una estructura algebraica que se construye como cociente de un producto directo $P$ de una familia de grupos $(G_i)_{i \in I}$. Donde el cociente est\'a dado por una cierta relaci\'on en $P$ que se define usando un filtro $\F$ en $I$.  Cuando el filtro $\F$ es ultrafiltro, a dicha estructura se le conoce como ultraproducto. Adem\'as, vamos a presentar una clase de grupos, en donde cada grupo de dicha clase tiene la propiedad de poderse encajar en un cierto ultraproducto. Esta secci\'on est\'a basado principalmente en la monograf\'ia de Ceccherini-Silberstein y Coornaert \cite[secci\'on 7]{b1}.\\

Sea $(G_i)_{i\in I}$ una familia de grupos y $\F$ un filtro en el conjunto de \'indices $I$. Sea $P$ el producto directo de la familia $(G_i)_{i\in I}$, es decir 
\[
P = \prod_{i\in I} G_i.
\]
 
 Sean ${\bf g} =(g_i)_{i\in I}$ y ${\bf h} =(h_i)_{i\in I}$ elementos de $P$. Decimos que ${\bf g} \sim_{\F} {\bf h}$ si el conjunto $\{i \in I\colon g_i = h_i \}$ pertenece a $\F$.
 Esta relaci\'on es de equivalencia como se muestra en los puntos (1)-(3) de la siguiente proposici\'on.
\begin{prop}\label{721}
Para todo ${\bf a}, {\bf b}, {\bf c}, {\bf d} \in P$ se tiene que:
\begin{enumerate}
\item ${\bf a} \sim_{\F} {\bf a}$;
\item ${\bf a} \sim_{\F} {\bf b}$ si y solo si ${\bf b} \sim_{\F} {\bf a}$;
\item Si ${\bf a} \sim_{\F} {\bf b}$ y ${\bf b} \sim_{\F} {\bf c}$ entonces ${\bf a} \sim_{\F} {\bf c}$;
\item Si ${\bf a} \sim_{\F} {\bf b}$ y ${\bf c} \sim_{\F} {\bf d}$ entonces ${\bf a}{\bf c} \sim_{\F} {\bf b}{\bf d}$;
\item ${\bf a} \sim_{\F} {\bf b}$ si y solo si ${\bf a}^{-1} \sim_{\F} {\bf b}^{-1}$.
\end{enumerate}
\end{prop}
Con los puntos (3) y (4) de la proposici\'on anterior se puede demostrar la siguiente
\begin{prop}\label{NF-normal} Sea $N_\F= \{ {\bf g} \in P \colon {\bf g} \sim_{\F} {\bf 1} \}$, 
en donde ${\bf 1}=(1_{G_i})_{i \in I}$, con $1_G$ la identidad del grupo $G$. Entonces $N_\F$ es un subgrupo normal de $P$.
\end{prop}

Como $N_{\F}$ es un grupo normal de $P$, el grupo cociente $P_{\F} = P / N_{\F}$ est\'a bien definido. Este grupo se llama  {\it producto reducido} de la familia de grupos $(G_i)_{i\in I}$ con respecto al filtro $\F$. Si $\F$ es un ultrafiltro, decimos que $P_{\F}$ es el {\it ultraproducto} de la familia de grupos $(G_i)_{i\in I}$ con respecto al ultrafiltro $\F$.
La siguiente proposici\'on implica que el conjunto cociente $P/_\sim{_\F}$ y el conjunto $P/N_\F$ son iguales.
\begin{prop}
\label{relclas}
Si ${\bf g}, {\bf h} \in P$, entonces 
\[{\bf g} N_{\F} = {\bf h} N_{\F} \Longleftrightarrow {\bf g} \sim_{\F} {\bf h}.
\]
 \end{prop}

Cuando $\F$ es un filtro principal, el grupo $P/N_{\F}$ no es interesante en el sentido de que es isomorfo a un grupo en la familia $(G_i)_{i\in I}$ como se muestra a continuaci\'on. Sea $\F_k$ el filtro principal en $I$ generado por $k$, entonces ${\bf g} \sim_{\F_k} {\bf 1}$ si y solo si $g_{k}=1_{G_{k}}$. Por lo tanto $N_{\F_k} = \{ {\bf g} \in P \colon g_k=1_{G_k} \}$. Si definimos la funci\'on $\phi \colon P \to G_k$ como $\phi({\bf g})=g_k$, tenemos que $\phi$ es un homomorfismo sobreyectivo  con $\ker(\phi)=N_{\F_k}$. Por el primer teorema de homomorfismo de grupos obtenemos que $P/N_{\F_k} \simeq G_k$. 
\subsection{Grupos encajables localmente}
En esta secci\'on vamos a definir a los grupos encajables localmente como primer ejemplo de grupos que se pueden aproximar, en un cierto sentido, por una clase de grupos. Primero algunas definiciones.

\begin{defi}
Una colecci\'on de grupos  $\mathscr{G}$ es una {\it clase de grupos} si satisface lo siguiente: si $G \in \mathscr{G}$ y $G'$ es un grupo isomorfo a $G$ entonces $G' \in \mathscr{G}$.
\end{defi}
Por ejemplo $\mathscr{G}$ puede ser la clase de grupos finitos, la clase de grupos nilpotentes, la clase de grupos solubles, etc. 

\begin{defi}
Sean $G$ y $H$ dos grupos. Dado un subconjunto finito $F \subseteq G$, una funci\'on $\phi \colon G \to H$ se llama {\it $F$-casi-homomorfismo} de $G$ en $H$ si satisface:
\begin{enumerate}
\item[i)] $\phi(f_1f_2)= \phi(f_1)\phi(f_2)$, para todo $f_1,f_2 \in F$;
\item[ii)] $\phi|_F$ es inyectiva.
\end{enumerate}
\end{defi}
\begin{defi}
Sea $\mathscr{G}$ una clase de grupos. Se dice que un grupo $\Gamma$ es {\it localmente encajable} en $\mathscr{G}$ si para todo subconjunto finito $F \subseteq \Gamma$, existe un grupo $G \in \mathscr{G}$ y un $F$-casi-homomorfismo $\varphi$ de $\Gamma$ en $G$.
\end{defi}

Por ejemplo, el grupo $\mathbb{Z}$ es localmente encajable en la clase de grupos finitos (grupos LEF por sus siglas en ingl\'es) como se muestra a continuaci\'on. Sea $F$ un subconjunto finito de $ \mathbb{Z}$, elegimos a un entero $n \geq 0$ tal que $F \subset [-n, n]$. Entonces, el homomorfismo cociente $\phi \colon \mathbb{Z} \rightarrow \mathbb{Z}/ (2n+1)\mathbb{Z}$ es un $F$-casi-homomorfismo. M\'as adelante veremos otros ejemplos de grupos LEF.

Algunas propiedades de cerradura de esta clase de grupos son las siguientes.
\begin{prop}\label{le-subgrupo}
Sea $\mathscr{G}$ una clase de grupos. Todo subgrupo de un grupo que es localmente encajable en $\mathscr{G}$ es localmente encajable en $\mathscr{G}$.
\end{prop}

\begin{prop}
Sea $\mathscr{G}$ una clase de grupos. Entonces $\Gamma$ es localmente encajable en $\mathscr{G}$ si y solo si todo subgrupo finitamente generado de $\Gamma$ es localmente encajable en $\mathscr{G}$.
\end{prop}

En el caso en que la clase $\mathscr{G}$ sea cerrada bajo productos directos tenemos que la encajabilidad local es cerrada al tomar producto directo. 
\begin{prop}
Sea $\mathscr{G}$ una clase de grupos la cual es cerrada bajo productos directos finitos. Sea $\left( \Gamma_i\right)_{i \in I}$ una familia de grupos que son localmente encajables en $\mathscr{G}$. Entonces, su producto directo $\Gamma = \prod_{i \in I} \Gamma_i$ es localmente encajable en $\mathscr{G}$.
\end{prop}

Los grupos localmente encajables son un ejemplo de grupos que se pueden definir usando ultraproductos como lo muestra el siguiente teorema.
\begin{thm}\label{gle-ultra-20}
Sea $\mathscr{G}$ una clase de grupos y sea $\Gamma$ un grupo. Las siguientes condiciones son equivalentes.
\begin{enumerate}
\item[a)] $\Gamma$ es localmente encajable en $\mathscr{G}$;
\item[b)] Existe una familia de grupos $(G_i)_{i \in I}$ tal que $G_i \in \mathscr{G}$, para todo $i \in I$, y un ultrafiltro $\U$ en $I$ tal que $\Gamma$ es isomorfo a un subgrupo del ultraproducto $P_{\U}$ de la familia $(G_i)_{i \in I}$ con respecto al ultrafiltro $\U$.
\end{enumerate}
\end{thm}

Si $\Gamma$ es un grupo que cumple con la condici\'on (b) del teorema anterior decimos que $\Gamma$ es {\it aproximado} por los grupos en la familia $(G_i)_{i \in I}$. Una consecuencia de este teorema es que la definici\'on de grupos encajables localmente no depende de la selecci\'on particular del ultrafiltro. 

Al parecer, Gordon y  Vershik fueron los primeros en estudiar a los grupos encajables localmente en la clase de grupos finitos (LEF) \cite{VeG}. En su trabajo ellos mostraron que las siguientes clases de grupos son LEF: grupos finitos, grupos libres, grupos abelianos, grupos nilpotentes, grupos de matrices. Ellos tambi\'en demostraron que 
cualquier grupo finitamente presentado con problema de la palabra no decidible (Pyotr Novikov~\cite{novi} demostr\'o que tales grupos existen) no es un grupo LEF.  Adem\'as, demostraron que el grupo
\[
G=\langle b, t\colon t^{-1}b^2t=b^3\rangle
\]
es un grupo finitamente presentado con problema de la palabra decidible que no es LEF.  En el libro de Magnus, Karrass y Solitar \cite{magnus} se pueden consultar las definiciones y resultados b\'asicos sobre el tema de presentaciones de grupos y sobre el problema de la palabra. 

\section{Aproximaci\'on m\'etrica de grupos}\label{sec-aproximacion-grupos}
En esta secci\'on vamos a presentar otro tipo de aproximaci\'on de grupos conocida como aproximaci\'on m\'etrica. Esta secci\'on est\'a basada principalmente en los art\'iculos de Stolz \cite[secci\'on 2]{artstolz} y de Stolz y Thom \cite[secci\'on 2]{stolz-thom}. Ambos art\'iculos forman parte de la tesis doctoral de Stolz \cite{tesisstolz}. Vamos a utilizar  funciones de longitud, en lugar de m\'etricas, que es en la manera como lo han hecho recientemente diversos autores \cite{lev-gleb, holt-rees, tesisstolz, artstolz, thom2, stolz-thom}. 
 
\begin{defi}
Sea $G$ un grupo. Una funci\'on $\ell \colon G \to [0, 1]$ es una funci\'on de longitud en $G$ si para todo $g, h \in G$ se cumple lo siguiente:

\begin{enumerate}
\item[FL1)] $\ell(g) = 0 \; \text{si y solo si} \; g=1$;
\item[FL2)] $\ell(g) = \ell(g^{-1})$;
\item[FL3)] $\ell(gh) \leq \ell(g) + \ell(h)$.
\end{enumerate}
\end{defi}

Decimos que $\ell$ es {\it invariante} si adem\'as se cumple $\ell(hgh^{-1})=\ell(g)$, para todo $h, g \in G$. Notemos que esta condici\'on es equivalente a pedir que $\ell(gh)=\ell(hg)$ para todo $h, g \in G$. A una funci\'on $\ell$ se le llama {\it pseudo funci\'on de longitud} si remplazamos la condici\'on (FL1) por $\ell(1)=0$.

La {\it funci\'on de longitud trivial} $\ell_0$ sobre un grupo $G$ se define como $\ell_0(g)=1$ si $g \neq 1$ y $\ell_0(1)=0$. Las (pseudo) funciones de longitud est\'an relacionadas con las (pseudo) m\'etricas como sigue  

\begin{prop}
Sea $\ell$ una (pseudo) funci\'on de longitud en un grupo $G$. Entonces $d(g,h):= \ell(gh^{-1})$ define una (pseudo) m\'etrica en $G$. Si $\ell$ es invariante entonces $d$ es bi-invariante. Si por el contrario $d$ es una (pseudo) m\'etrica izquierdo-invariante en $G$ entonces $\ell(g):= d(g,1)$ es una (pseudo) funci\'on de longitud. Si $d$ es bi-invariante entonces $\ell$ es invariante.
\end{prop}

Algunos ejemplos de funciones de longitud son:
\begin{enumerate}
\item 
Sea [n] el conjunto de los naturales $\{1, \dots, n\}$. Denotemos por  $S_n$ al grupo sim\'etrico $\Sym([n])$ que consiste de todas las funciones biyectivas de $[n]$ sobre $[n]$ con la operaci\'on de grupo dada por la composici\'on de funciones. La funci\'on de longitud de Hamming de una permutaci\'on $\sigma \in S_n$, se define como
\[
\ell_H(\sigma)= \frac{|\{i \in [n] \colon \sigma(i)\neq i\}|}{n}.
\]
\item 
El grupo $U_n$ de elementos unitarios de $M_n(\CC)$ se puede equipar con la funci\'on de longitud invariante de Hilbert-Schmidt definida como
\[
\ell_{HS}(g)=d_{HS}(g, 1),
\]
con $d_{HS}(A, B)= \sqrt{\tau\left((A-B)^*(A-B)\right)}$, en donde $\tau$ es la traza normalizada en $U_n(\CC)$ dada por $\tau(C)= \frac{1}{n}\sum\limits_{i=1}^n c_{ij}$, y $C^*$ es la matriz transpuesta conjugada de $C$, para $C \in U_n(\CC)$.
\item 
En $GL(n,\mathbb{C})$, el grupo de las matrices de $n\times n$ invertibles sobre $ \mathbb{C}$, se define la funci\'on de longitud del rango
\[
\ell_{L}(A)=\frac{1}{n}rango(I-A).
\]

\item 
Para cualquier grupo finito $G$ se puede definir la pseudo funci\'on de longitud conjugaci\'on dada por
\[
\ell_c(g) = \frac{log|C(g)|}{log|G|}
\]
donde $C(g)$ denota la clase de conjugaci\'on de $g$.
\end{enumerate}
La siguiente proposici\'on muestra una manera de definir una pseudo funci\'on de longitud invariante en el cociente de un grupo que tiene pseudo funci\'on de longitud invariante. 
\begin{prop}
Sea $G$ un grupo con pseudo funci\'on de longitud invariante $\ell$ y $H$ un subgrupo normal de $G$. Entonces 
\[\ell_{G/H}(gH) := \inf_{x \in H} \ell(gx)
\] define una pseudo funci\'on de longitud invariante en $G/H$. Si $G$ es finito y $\ell$ es una funci\'on de longitud, entonces $\ell_{G/H}$ es una funci\'on de longitud.
\end{prop}

\subsection{Ultraproductos m\'etricos de grupos}
Ahora vamos a construir los ultraproductos m\'etricos. Para ello necesitamos el siguiente resultado. 
\begin{prop}\label{n-normal}
Sea $G$ un grupo con pseudo funci\'on de longitud invariante $\ell$. Entonces el conjunto 
\[ 
N= \{g \in G \colon \ell(g)=0 \}
\] es un subgrupo normal de $G$. Adem\'as 
\[
\ell_{G/N}(gN)=\ell(g), \; g \in G,
\] y $\ell_{G/N}$ define una funci\'on de longitud invariante en $G/N$.
\end{prop}
Para que $N$ sea un grupo normal es necesario que la pseudo funci\'on de longitud sea invariante como lo muestra el siguiente ejemplo, que aparece en el resumen de Pestov~\cite{pes}. Para este ejemplo vamos a utilizar m\'etrica en lugar de funci\'on de longitud. 
 
\begin{ejem}
Sea $S_{\mathbb{N}}$ el grupo sim\'etrico que consiste de todas las biyecciones de $\mathbb{N}$ sobre si mismo, equipado con la siguiente m\'etrica:

\[ 
d(\sigma, \tau)= \sum_{i \in D} 2^{-i}
\]
en donde $D=\{ i \in \NN \colon \sigma(i) \neq \tau(i)\}$, y usamos la convenci\'on de que la suma es cero cuando $D$ es el conjunto vac\'io. Sea $\U$ un ultrafiltro libre sobre $\mathbb{N}$. Sean $\sigma, \tau \in \prod_\mathbb{N} S_{\mathbb{N}}$, tales que $\sigma = \left( \sigma_i \right)$, con $\sigma_i=  \left( i, i+1 \right)$, y $\tau = \left( \tau_i \right)$, con $\tau_i= \left( 1, i \right)$ (las funciones $\sigma_i$ y $\tau_i$ est\'an expresadas en notaci\'on c\'iclica). Es decir, $\sigma$ y $\tau$ son dos sucesiones de transposiciones en $S_{\mathbb{N}}$. Puesto que $d\left( \sigma_i, 1 \right)= 2^{-i}+2^{-(i+1)}$ tenemos que $ \lim\limits_{u} d\left( \sigma_i, 1 \right)= 0$ y as\'i $\sigma \in N$. Pero, como $\tau_i \sigma_i\tau_i^{-1}= \left(1, i+1 \right)$ y $d\left( \left(1, i+1 \right), 1 \right)= 2^{-1}+2^{-(i+1)}$, entonces $\lim\limits_u d\left( \sigma_i^{-1}\tau_i\sigma_i, 1 \right)= 1/2$, lo cual implica que $\tau_i \sigma_i\tau_i^{-1} \not \in N$. Por lo tanto $N$ no es un subgrupo normal de $\prod_\mathbb{N} S_{\mathbb{N}}$.
\end{ejem}

En la siguiente proposici\'on se define una pseudo funci\'on de longitud invariante en el producto directo de grupos equipados con pseudo funci\'on de longitud invariante.

\begin{prop}\label{prod-direc-lon}
Sean $(G_i,\ell_i)_{i\in I}$ una secuencia de grupos cada uno equipado con una pseudo funci\'on de longitud $\ell_i$ invariante. Sea  $P = \prod_{i\in I} G_i$ su producto directo. Sea $\U$ un ultrafiltro libre sobre $I$. Para cada ${\bf g} \in P$, definimos \[
\tilde{\ell}({\bf g}):= \lim_{u} \ell_i(g_i).\] Entonces $\tilde{\ell}$ define una pseudo funci\'on de longitud invariante en $P$.
 \end{prop}

Por las proposiciones \ref{n-normal} y \ref{prod-direc-lon} obtenemos que el grupo cociente 
\[
(P)_u:=\prod_{i\in I} G_i/N,
\]
est\'a bien definido, con $N=\{{\bf g}\in P \colon \tilde{\ell}({\bf g})=0\}$, y adem\'as $\ell_{P/N}({\bf g}N)$ es una funci\'on de longitud invariante en $P/N$ tal que $\ell_{P/N}({\bf g}N) = \tilde{\ell}({\bf g})$, para todo ${\bf g} \in P$. 

Al grupo $(P)_u$ se lo conoce como el {\it ultraproductro m\'etrico} de la familia $(G_i,\ell_i)_{i\in I}$  m\'odulo $\U$.

\subsection{Grupos $\mathscr{G}$--aproximables}\label{grupos-g-aprox}

Vamos a dar una definici\'on de aproximaci\'on m\'etrica de grupos que usa ultraproductos.
\begin{defi}\label{g-aprox-ultra}
Sea $\mathscr{G}$ una clase de grupos, en donde cada grupo en la clase $\mathscr{G}$ est\'a equipado con una pseudo funci\'on de longitud invariante. Decimos que un grupo $\Gamma$ tiene la propiedad de $\mathscr{G}$-aproximaci\'on m\'etrica si existe un conjunto de \'indices $I$ y un ultrafiltro $\U$ en $I$ tal que $\Gamma$ es isomorfo a un subgrupo de un ultraproducto m\'etrico $(\prod_{i \in I} G_i)_u$, con grupos $G_i \in \mathscr{G}$. En este caso decimos que el grupo $\Gamma$ es aproximado de forma m\'etrica por los grupos en la clase $ \mathscr{G}$.
\end{defi}
Desde hace varios a\~nos se han estudiado clases de grupos que tienen la propiedad de aproximaci\'on m\'etrica. A estos grupos se les han dado diferentes nombres dependiendo de la clase de grupos $\mathscr{G}$ que los apro\-xi\-man y de las m\'etricas con las que cuentan los grupos en $\mathscr{G}$.
\begin{defi}\label{muchas-defi}
Si $\Gamma$ tiene la propiedad de $\mathscr{G}$-aproximaci\'on, al grupo $\Gamma$ se le llama:
\begin{itemize}
\item Localmente encajable en finito si $\mathscr{G}$ es la clase de grupos finitos con la funci\'on de longitud trivial $\ell_0$ (Gordon y Vershik, \cite{VeG}). 

\item  S\'ofico si $\mathscr{G}$ es la clase de grupos sim\'etricos finitos con la funci\'on de longitud de Hamming (Gromov, \cite{gromov-artc}).  
\item Hiperlineal si $\mathscr{G}$ es la clase de matrices unitarias de rango finito sobre $\CC$ con la funci\'on de longitud de Hilbert-Schmidt (R\u{a}dulescu, \cite{radu}).
\item S\'ofico d\'ebil si $\mathscr{G}$ es la clase de todos los grupos finitos, en donde cada grupo est\'a equipado con alguna funci\'on de longitud invariante (Glebsky y Rivera, \cite{docluis}).
\item S\'ofico lineal si $\mathscr{G}$ es la clase de grupos lineales generales $GL(n, \CC)$ con $\ell(A)=\frac{1}{n} rango(I-A)$ (Arzhantseva y P\u{a}unescu, \cite{arzpau}).
\item $K$-s\'ofico si $\mathscr{G}$ es la clase de grupos lineales generales $GL(n, F)$, para $F$ un campo fijo, con la longitud rango (Stolz, \cite{artstolz}).
\end{itemize}
\end{defi}

 Los grupos s\'oficos est\'an relacionados con los grupos hiperlineales de la siguiente manera. A cada permutaci\'on $\sigma \in S_n$ le podemos asociar una matriz $A_{\sigma}$ de $n \times n$ de la siguiente manera:
\[
(A_{\sigma})_{ij}= \left\{
\begin{array}{cl}
1& \text{si}\; \sigma(j)=i\\
0& \text{en otro caso}
\end{array}
\right.
\]

La funci\'on $\sigma \mapsto A_\sigma$ define un encaje de $S_n$ al grupo unitario $U(n)$. Se puede verificar que $\ell_H(\alpha)=\frac{1}{2}\left(\ell_{HS}(A_{\alpha})\right)^2$, para toda $\alpha$ en $S_n$. Usando lo anterior Elek y Szab\'o~\cite{ES2} demostraron que los grupos s\'oficos son hiperlineales. Las otras clases se relacionan como sigue: los grupos LEF son s\'oficos, los grupos s\'oficos son s\'oficos d\'ebiles  \cite{docluis}, los grupos s\'oficos son s\'oficos lineales, y los grupos s\'oficos lineales son s\'oficos d\'ebiles  \cite{arzpau}. En todos los casos se desconoce si el rec\'iproco es verdadero. A la fecha no se conocen ejemplos de grupos que no sean s\'oficos, o hiperlineales, o s\'oficos d\'ebiles, o s\'oficos lineales, o $K$-s\'oficos, y el encontrar un ejemplo de un grupo que no cumpla con alguna de estas definiciones es uno de los problemas m\'as importantes en la teor\'ia de aproximaci\'on de grupos. En la siguiente secci\'on continuaremos con la exposici\'on sobre grupos s\'oficos.

\section{Grupos s\'oficos}\label{grupos-soficos}
Los grupos s\'oficos, son una de las clases de grupos $\mathscr{G}$-aproximables m\'as estudiadas a la fecha y son una generalizaci\'on com\'un de dos clases importantes de grupos: los grupos residualmente finitos y los grupos amenables. En esta secci\'on vamos a dar otra definici\'on de grupos s\'oficos que no usa ultraproductos y vamos a presentar algunas propiedades de cerradura de esta clase de grupos. La equivalencia de las dos definiciones la da el teorema \ref{equivdef}. 

\begin{defi}
Sea $G$ un grupo, $K \subseteq G$ un subconjunto finito no vac\'io, $\epsilon > 0$ y $F$ un conjunto finito no vac\'io. Una funci\'on $\varphi: G \rightarrow \Sym(F)$ se llama {\it $(K, \epsilon)$-casi-homomorfismo} si satisface las siguientes condiciones:
\begin{enumerate}
\item[i)] para todo $g, h \in K$, $d_H(\varphi(g)\varphi(h), \varphi(gh))\leq \epsilon$;
\item[ii)] para todo $g, h \in K$, $g \neq h$, se tiene $d_H(\varphi(g), \varphi(h))\geq 1- \epsilon$.
\end{enumerate}
\end{defi}
\begin{defi}

Un grupo $\Gamma$ es llamado s\'ofico si para todo $K \subseteq \Gamma$ subconjunto finito y para todo $\epsilon > 0$, existe un conjunto finito no vac\'io $F$ y un $(K, \epsilon)$-casi-homomorfismo $\varphi: G \rightarrow \Sym(F)$.
\end{defi}
Para el caso de grupos finitamente generados esta definici\'on de grupos s\'oficos es equivalente a la definici\'on original presentada por Gromov \cite{gromov-artc} y Weiss \cite{weiss} (una demostraci\'on se puede consultar, por ejemplo, en~\cite[secci\'on 7.7]{b1}). La siguiente proposici\'on aparece en \cite[cap\'itulo 6]{connesQWEP} y es una recopilaci\'on de otras definiciones equivalentes de grupos s\'oficos que han aparecido a lo largo de la literatura. 

\begin{prop}\label{equiv-sofico}
Sea $\Gamma$ un grupo. Entonces las siguientes condiciones son equivalentes:
\begin{enumerate}
\item Para cada subconjunto finito $K \subseteq \Gamma$ y para todo $\epsilon >0$,  existe un conjunto finito no vac\'io $F$ y una funci\'on $\phi \colon \Gamma \rightarrow \Sym(F)$ con las propiedades:
\begin{enumerate}
\item[i)] $d_H(\phi(g)\phi(h), \phi(gh)) \leq \epsilon$ para cada $g, h \in K$;
\item[ii)] $d_H(\phi(1_{\Gamma}), 1_{\Sym(F)}) \leq \epsilon$;
\item [iii)]  $d_H(\phi(g),\phi(h)) \geq 1- \epsilon$ para cada $g, h \in K$ con $g \neq h$.\\
(El grupo $\Gamma$ tiene la propiedad de $\mathscr{G}$-aproximaci\'on discreta fuerte)
\end{enumerate}
\item Para cada subconjunto finito $K \subseteq \Gamma$ y para todo $\epsilon >0$,  existe un conjunto finito no vac\'io $F$ y una funci\'on $\phi \colon \Gamma \rightarrow \Sym(F)$ con las propiedades:
\begin{enumerate}
\item[i)] $d_H(\phi(g)\phi(h), \phi(gh)) \leq \epsilon$ para cada $g, h \in K$;
\item [ii)]  $d_H(\phi(g),\phi(h)) \geq 1- \epsilon$ para cada $g, h \in K$ con $g \neq h$.
\end{enumerate}
\item Para cada constante $\delta \in (0,1)$, cada subconjunto finito $K \subseteq \Gamma$ y todo $\epsilon >0$, existe un conjunto finito no vac\'io $F$ y una funci\'on $\phi \colon \Gamma \rightarrow \Sym(F)$ con las propiedades:
\begin{enumerate}
\item[i)] $d_H(\phi(g)\phi(h), \phi(gh)) \leq \epsilon$ para cada $g, h \in K$;
\item [ii)]  $d_H(\phi(g),\phi(h)) \geq \delta$ para cada $g, h \in K$ con $g \neq h$.\\
(El grupo $\Gamma$ tiene la propiedad de $\mathscr{G}$-aproximaci\'on discreta)
\end{enumerate}

\item Existe alg\'un $\delta >0$ tal que para cada subconjunto finito $K \subseteq \Gamma$ y todo $\epsilon >0$, existe un conjunto finito no vac\'io $F$ y una funci\'on $\phi \colon \Gamma \rightarrow \Sym(F)$ con las propiedades:
\begin{enumerate}
\item[i)] $d_H(\phi(g)\phi(h), \phi(gh)) \leq \epsilon$ para cada $g, h \in K$;
\item [ii)]  $d_H(\phi(g),\phi(h)) \geq \delta$ para cada $g, h \in K$ con $g \neq h$.\\
($\Gamma$ tiene la propiedad de $\mathscr{G}$-aproximaci\'on)
\end{enumerate}
\item Para cada conjunto finito $K \subseteq \Gamma$, existe un $\delta_K$ tal que para cada $\epsilon >0$, existe un subconjunto finito no vac\'io $F$ y una funci\'on $\varphi \colon \Gamma \to \Sym(F)$ con las siguientes propiedades:
\begin{enumerate}
\item [i)] $d_H(\phi(g)\phi(h), \phi(gh)) \leq \epsilon$ para cada $g, h \in K$;
\item [ii)]$d_H(\phi(g),\phi(h)) \geq \delta_K$ para cada $g, h \in K$ con $g \neq h$.
\end{enumerate}
\end{enumerate}
\end{prop}
\begin{defi}
Un grupo $\Gamma$ se llama {\it s\'ofico} si cumple alguna (y por lo tanto todas) de las condiciones (1)-(5) de la proposici\'on anterior. 
\end{defi}
\begin{thm}
La clase de grupos s\'oficos es cerrada con respecto a las siguientes operaciones:
\begin{enumerate}
\item Subgrupos \cite{ESprodsof};
\item Producto directo \cite{ESprodsof};
\item L\'imites directos \cite{ESprodsof};
\item L\'imites inversos \cite{ESprodsof};
\item Productos libres \cite{ESprodsof};
\item Producto gr\'afico \cite{laura};
\item Producto trenzado \cite{hayes};
\item Producto libre amalgamado sobre grupos amenables \cite{ES31, paunescu} ;
\item Extensiones por grupos amenables: si $N \lhd G$, $N$ es s\'ofico y $G/N$ es amenable, entonces $G$ es s\'ofico \cite{ESprodsof}; 
\item Extensiones HNN sobre grupos amenables  \cite{collis, ES31, paunescu}.
\end{enumerate}
\end{thm}

\section{Ejemplos de grupos s\'oficos}

Los grupos s\'oficos generalizan a dos clases importantes de grupos: los grupos residualmente finitos y los grupos amenables.

\begin{defi} 
Un grupo $\Gamma$ es {\it residualmente finito} si para cada elemento $g \in \Gamma$ con $g \neq 1_{\Gamma}$, existe un grupo finito $G$ y un homomorfismo $\phi \colon \Gamma \rightarrow G$ tal que $\phi(g) \neq 1_G$.
\end{defi}

La clase de grupos residualmente finitos contienen a todos los grupos finitos, a los grupos abelianos finitamente generados (en particular el grupo aditivo $\mathbb{Z}$), y a los grupos libres, entre otros. La clase de grupos residualmente finitos es cerrada al tomar subgrupos y l\'imites inversos. Un teorema de Mal'cev \cite{aimal1} muestra que todos los grupos lineales finitamente generados son residualmente finitos. Los grupos aditivos $\Q$, $\RR$ y $\CC$ no son grupos residualmente finitos. En el libro \cite[secci\'on 2]{b1}  se puede consultar las demostraciones de lo dicho anteriormente, y otras propiedades de esta clase de grupos. Dos problemas importantes que permanecen abiertos sobre este tema son los siguientes: 1) determinar si todos los grupos Gromov-hiperb\'olicos son residualmente finitos; 2) determinar bajo que condiciones un grupo con un solo relator es re\-si\-dual\-men\-te finito \cite{sapir}.  

\begin{prop}
Todo grupo residualmente finito es s\'ofico.
\end{prop}

\begin{defi} 
Un grupo discreto $\Gamma$ es amenable si para todo subconjunto finito $\Lambda \subset \Gamma$  y para todo $\epsilon >0$ existe un conjunto $E$ que es $(\Lambda, \epsilon)$-Følner, es decir, un subconjunto finito $E$ de $\Gamma$ tal que para todo $g \in \Lambda$ 
\[ | gE \bigtriangleup E| < \epsilon|E|,\]
 donde $\bigtriangleup$ denota la diferencia sim\'etrica entre conjuntos. 
\end{defi}

Cualquier grupo finito es amenable: sea $\Gamma$ un grupo finito, entonces para cualquier subconjunto finito $\Lambda \subset \Gamma$  y para todo $\epsilon >0$, el grupo $\Gamma$ es un conjunto $(\Lambda, \epsilon)$-Følner porque $|g\Gamma \Delta \Gamma|=0$.

La clase de los grupos amenables contiene a los grupos finitos, a los grupos abelianos y a los grupos solubles, entre otros. El grupo $\Q$ es amenable pero no residualmente finito y el grupo libre con dos ge\-ne\-ra\-do\-res $F_2$ es residualmente finito pero no amenable (las demostraciones de estos hechos y m\'as informaci\'on de esta clase de grupos se puede consultar en \cite[secci\'on 4]{b1}). 

\begin{prop}
Todo grupo amenable es s\'ofico.
\end{prop}

Finalizamos con el siguiente resultado que tiene como consecuencia que los grupos encajables localmente en la clase de grupos finitos (LEF) y los grupos encajables localmente en la clase de grupos amenables (LEA) son s\'oficos. 
\begin{prop}\label{loc-enc-sof}
Todo grupo $G$ que es localmente encajable en la clase de grupos s\'oficos es s\'ofico. \end{prop}

\section{Propiedades de grupos $\mathscr{G}$-aproximables}\label{2-prop-grupos-aprox}

En la secci\'on anterior se presenta una definici\'on de grupos s\'oficos que no utiliza ultrafiltros. Es una pregunta natural si es o no posible tener otra definici\'on de grupos $\mathscr{G}$-aproximables que no use ultraproductos. En esta secci\'on vamos a presentar dicha definici\'on y veremos que es muy similar a la de grupos s\'oficos. Esta secci\'on est\'a basado principalmente en \cite{holt-rees, tesisstolz, artstolz}.  

\begin{defi}\label{G-aprox-delta}
Sea $\mathscr{G}$ una clase de grupos, cada uno de los cuales est\'a equipado con una pseudo funci\'on de longitud $\ell$. Un grupo $\Gamma$ tiene la propiedad de $\mathscr{G}$-aproximaci\'on si para todo $g \in \Gamma$, $g \neq 1$, existe $\delta_g > 0$ tal que para todo $\epsilon > 0$ y cualquier subconjunto finito $F \subseteq \Gamma$, $F \neq \emptyset$, existe un grupo $G \in \mathscr{G}$ y una funci\'on $\varphi \colon \Gamma \to G$ tal que:

\begin{enumerate}
\item [{\it GA1})] $\ell(\varphi(1)) \leq \epsilon$;
\item [{\it GA2})] $\ell(\varphi(g)) \geq \delta_g$, para todo $g \in F \setminus \{1\}$;
\item [{\it GA3})] $\ell(\varphi(g)\varphi(h)\varphi(gh)^{-1}) \leq \epsilon$, para todo $g, h \in F$.
\end{enumerate}
A la funci\'on $\varphi$ se le conoce como un $(F, \epsilon, \delta_g, \ell)$-casi-homomorfismo. Cuando no es necesario hacer referencia a $\delta_g$ y a la pseudo funci\'on de longitud $\ell$ se dice simplemente que $\varphi$ es un $(F, \epsilon)$-casi-homomorfismo. 
\end{defi}

Se asume que $\delta_g \leq 1$, para todo $g \in \Gamma$ porque de lo contrario no se cumplir\'ia la condici\'on (GA2). 

Un grupo que tiene la propiedad de $\mathscr{G}$-aproximaci\'on se llama s\'ofico, hiperlineal, s\'ofico d\'ebil, s\'ofico lineal, $K$-s\'ofico dependiendo de quien es la clase $\mathscr{G}$ y la funci\'on de longitud de manera similar como en la definici\'on \ref{muchas-defi}.

Una funci\'on $\delta \colon \Gamma \to \RR$ tal que $\delta(g) > 0$, para $g \neq 1$ y $\delta(1)=1$ se le llama {\it funci\'on de peso} para el grupo $\Gamma$.  Notemos que en la definici\'on de $\mathscr{G}$-aproximaci\'on se requiere que el grupo $\Gamma$ tenga una funci\'on de peso.

Otros tipos de aproximaci\'on de grupos son, la $\mathscr{G}$-aproximaci\'on dis\-cre\-ta la cual se da al reemplazar en la definici\'on $\delta_g$ por una constante $\delta$; y la $\mathscr{G}$-aproximaci\'on fuerte donde en la condici\'on (GA2) se tiene $\ell(\varphi(g)) \geq diam(G) - \epsilon$ para todo $g \in F \setminus \{1\}$, en donde $diam(G)=\sup\limits_{g \in G} \ell(g)$, y a $\varphi$  se le conoce como {\it casi-homomorfismo fuerte}.

No se sabe si las diferentes versiones de $\mathscr{G}$-aproximaci\'on son e\-qui\-va\-len\-tes de manera general. En el caso de grupos s\'oficos e hiperlineales dichas versiones si son equivalentes. Para el caso de grupos s\'oficos la equivalencia la da la proposici\'on~\ref{equiv-sofico} y para el caso de grupos hiperlineales la equivalencia la da un resultado que se puede encontrar, por ejemplo, en la tesis de Olesen \cite[proposici\'on 5.1.2]{connesQWEP}. 

El siguiente resultado, bien conocido en el \'area, muestra la equivalencia de las definiciones de aproximaci\'on \ref{g-aprox-ultra} y \ref{G-aprox-delta} lo cual implica que la definici\'on de $\mathscr{G}$-aproximaci\'on no depende de la selecci\'on del ultrafiltro. 

\begin{thm}\label{equivdef}
Un grupo $\Gamma$ tiene la propiedad de $\mathscr{G}$-aproximaci\'on seg\'un la definici\'on~\ref{G-aprox-delta} si y solo si existe un conjunto de \'indices $I$ y un ultrafiltro $\U$ en $I$ tal que $\Gamma$ puede ser encajado en un ultraproducto $\left(\prod_{i \in I} G_i\right)_u$ de grupos $G_i \in \mathscr{G}$. 
\end{thm}

\subsection{Propiedades generales de grupos $\mathscr{G}$-aproximables}

Se conocen pocos resultados generales para los grupos $\mathscr{G}$-aproximables. Vamos a concluir este art\'iculo mencionando algunos de estos resultados.  Se sabe que la clase de grupos con la propiedad de  $\mathscr{G}$-aproximaci\'on es cerrada al tomar subgrupos y l\'imites inversos \cite[proposici\'on 3.3]{artstolz}. Lo mismo es cierto para los grupos con la propiedad discreta y la propiedad discreta fuerte de  $\mathscr{G}$-aproximaci\'on. Se sabe que si la clase $\mathscr{G}$ tiene una propiedad conocida como amplificaci\'on, entonces la clase de grupos $\mathscr{G}$-aproximables es cerrada al  tomar l\'imites directos \cite[proposici\'on 3.5]{artstolz}.


Se conocen pocos resultados para el caso del producto directo de grupos $\mathscr{G}$-aproximables. Stolz demostr\'o un resultado parcial con varias hip\'otesis \cite[proposici\'on 3.4]{artstolz}. Vamos a presentar el resultado que Derek F. Holt y Sarah Rees \cite{holt-rees}    anuncian en 2016. 

Primero vamos a definir una funci\'on de longitud en el producto directo de dos grupos con funci\'on de longitud. Si $\mathscr{G}$ es una clase de grupos m\'etricos, escribiremos $(G, \ell_G) \in \mathscr{G}$ para indicar que $G \in \mathscr{G}$ y que $\ell_G$ es una funci\'on de longitud definida en $G$. 
\begin{prop}
Sea $\mathscr{G}$ una clase de grupos m\'etricos. Supongamos que $(G, \ell_G), (H, \ell_H) \in \mathscr{G}$. Entonces para $p \in \mathbb{N} \cup \{\infty\}$ y $\alpha \in G \times H$ con ${\bf a}=(g, h)$, la funci\'on $\ell_{G \times H}^p \colon G \times H \rightarrow [0,1]$ definida por

\[
\ell_{G \times H}^p(g, h) = \sqrt[p]{\frac{\ell_G(g)^p +\ell_H(h)^p}{2}}, \qquad p \in \mathbb{N}
\]

y 

\[
\ell_{G \times H}^\infty (g, h) = \max \{ \ell_G(g), \ell_H(h) \}
\]

es una funci\'on de longitud invariante. 
\end{prop}

\begin{thm}\label{prod-claseG}
Sea $\mathscr{G}$ una clase de grupos con funciones de longitud invariantes asociadas y supongamos que, para alg\'un $p \in \mathbb{N} \cup \{ \infty\}$ fijo, y para cualquier par de grupos $G, H \in \mathscr{G} $ se cumple que
\[
\left(G, \ell_G \right), \left(H, \ell_H \right) \in \mathscr{G} \Rightarrow \left( G \times H, \ell^p \right) \in \mathscr{G}.
\]
Entonces el  producto directo $G \times H$ de dos grupos $G$ y $H$ $\mathscr{G}$-aproximables es $\mathscr{G}$-aproximable.
\end{thm}

El resultado de este teorema puede ser aplicado para deducir la ce\-rra\-du\-ra bajo productos directos para las clases de grupos s\'oficos d\'ebiles, grupos LEF, grupos hiperlineales, grupos lineales s\'oficos, donde la condici\'on se mantiene de la siguiente manera:
\begin{itemize}
\item Grupos s\'oficos d\'ebiles para toda $p$;
\item Grupos LEF para $p = \infty$;
\item Grupos hiperlineales para $p=2$;
\item Grupos lineales s\'oficos para $p=1$.
\end{itemize}

Para la funci\'on de longitud de Hamming  $\ell_G$, $\ell_H$, la funci\'on $\ell^p_{\ell_G, \ell_h}$ no es una funci\'on de longitud de Hamming, y por lo tanto no podemos deducir la cerradura de la clase de grupos s\'oficos bajo el producto directo a partir de este resultado.

Concluimos este trabajo citando los art\'iculos de Holt y Ress \cite{holt-rees}, y de Hayes y Sale \cite{hayes2} en donde se pueden encontrar otras propiedades de cerradura para los grupos $\mathscr{G}$-aproximables. 

\bigskip
{\bf \centerline{Agradecimientos}}
Este art\'iculo es parte de la tesis de licenciatura  de la segunda autora bajo la direcci\'on del primer autor. La tesis se present\'o en la Unidad Acad\'emica de Matem\'aticas de la Universidad Aut\'onoma de Zacatecas en marzo del 2016. Los autores agradecen a los profesores Daniel Duarte, Patricia Jim\'enez, Jes\'us Lea\~nos y Alexander Pyshchev por sus comentarios y sugerencias sobre la tesis, que se ven reflejados en este art\'iculo. Los autores agradecen al editor y al revisor por sus \'utiles sugerencias y correcciones.


\bigskip
\hfill\
{\footnotesize
\parbox{5.2cm}{Luis Manuel Rivera Mart\'inez\\
{\it Unidad Acad\'emica de Matem\'aticas},\\
Universidad Aut\'onoma de Zacatecas,\\
Zacatecas, M\'exico\\
{\sf luismanuel.rivera@gmail.com}}\
{\hfill}\
\parbox{5.2cm}{Nidya Montserrath Veyna Garc\'ia\\
{\it Unidad Acad\'emica de Matem\'aticas},\\
Universidad Aut\'onoma de Zacatecas,\\
Zacatecas, M\'exico\\
{\sf nidyaveyna@gmail.com }}}\
\hfill


\end{document}